\documentstyle[11pt,fullpage,twoside,amssym,saks]{article}
\def\Xint#1{\mathchoice
   {\XXint\displaystyle\textstyle{#1}}%
   {\XXint\textstyle\scriptstyle{#1}}%
   {\XXint\scriptstyle\scriptscriptstyle{#1}}%
   {\XXint\scriptscriptstyle\scriptscriptstyle{#1}}%
   \!\int}
\def\XXint#1#2#3{{\setbox0=\hbox{$#1{#2#3}{\int}$}
     \vcenter{\hbox{$#2#3$}}\kern-.5\wd0}}
\def\dashint{\Xint-}

\begin{document}
\markboth{\centerline{ RUSLAN SALIMOV}} {\centerline{ ON FINITELY
LIPSCHITZ    SPACE MAPPINGS.}}

\def\kohta #1 #2\par{\par\noindent\rlap{#1)}\hskip30pt
\hangindent30pt #2\par}
\def\esssup{\operatornamewithlimits{ess\,sup}}
\def\tomes{\mathop{\longrightarrow}\limits^{mes}}
\def\ts{\textstyle}
\def\I{\roman{Im}}
\def\mes{\mbox{\rm mes}}
\def\A{{{\cal {A}}}}
\def\Rm{{{\Bbb R}^m}}
\def\Rn{{{\Bbb R}^n}}
\def\Rk{{{\Bbb R}^k}}
\def\R3{{{\Bbb R}^3}}
\def\lR{{\overline {{\Bbb R}}}}
\def\lRn{{\overline {{\Bbb R}^n}}}
\def\lRm{{\overline {{\Bbb R}^m}}}
\def\lRk{{\overline {{\Bbb R}^k}}}
\def\lBn{{\overline {{\Bbb B}^n}}}
\def\Bn{{{\Bbb B}^n}}
\def\R{{\Bbb R}}
\def\Z{{\Bbb Z}}
\def\C{{\Bbb C}}
\def\B{{\Bbb B}}
\def\e{{\varepsilon}}
\def\L{{\Lambda}}
\def\f{{\varphi}}
\def\F{{\Phi}}
\def\x{{\chi}}
\def\d{{\delta }}
\def\D{{\Delta }}
\def\c{{\circ }}
\def\tg{{\tilde{\gamma}}}
\def\a{{\alpha }}
\def\p{{\psi }}
\def\m{{\mu }}
\def\r{{\rho }}
\def\S{{\Sigma }}
\def\O{{\Omega }}
\def\s{{\sigma }}
\def\l{{l}\, }
\def\g{{\gamma }}
\def\G{{\Gamma }}
\def\D{{\Delta }}
\let\text=\mbox
\let\Cal=\cal

\def\cc{\setcounter{equation}{0}
\setcounter{figure}{0}\setcounter{table}{0}}

\overfullrule=0pt

\def\eqb{\begin{equation}}
\def\eqe{\end{equation}}
\def\eb{\begin{eqnarray}}
\def\ee{\end{eqnarray}}
\def\ebnn{\begin{eqnarray*}}
\def\eenn{\end{eqnarray*}}
\def\db{\begin{displaystyle}}
\def\de{\end{displaystyle}}
\def\tb{\begin{textstyle}}
\def\te{\end{textstyle}}
\def\exb{\begin{ex}}
\def\exe{\end{ex}}
\def\bth{\begin{theo}}
\def\eth{\end{theo}}
\def\bcor{\begin{corol}}
\def\ecor{\end{corol}}
\def\blem{\begin{lemma}}
\def\elem{\end{lemma}}
\def\brem{\begin{rem}}
\def\erem{\end{rem}}
\def\bpr{\begin{propo}}
\def\epr{\end{propo}}

\title{{\bf On finitely  Lipschitz  space mappings.}}

\author{{\bf  R. Salimov}\\}
\date{\today \hskip 4mm ({\tt SALIMOV.tex})}
\maketitle

\large \abstract

It is established that a ring $Q$-homeomorphism with respect to
$p$-modulus  in ${\Bbb R}^n$, $n\geqslant2,$ is finitely Lipschitz
if $n-1<p<n$ and if the mean integral value of $Q(x)$ over
infinitisimial balls $B(x_0,\varepsilon)$ is finite everywhere.
\endabstract
\bigskip



\bigskip
\large \cc
\section{Introduction}

Recall that, given a family of paths $\Gamma$ in $\Bbb{R}^n$, a
Borel function $\varrho:\Bbb{R}^n\to[0,\infty]$ is called {\it
admissible} for $\Gamma$, abbr. $\varrho\in adm\,\Gamma$, if

\eqb \label{eq13.2} \int\limits_{\gamma}\varrho\,ds\ \geq\ 1
 \eqe for all
$\gamma\in\Gamma$. The  {\it $p$-modulus} of $\Gamma$ is the
quantity \eqb\label{eq13.3}M_{p}(\Gamma)\ =\ \inf\limits_{\varrho\in
adm\,\Gamma}\int\limits_{\Rn}\varrho^{p}(x)\ dm(x)\ . \eqe Here the
notation $m$ refers to the Lebesgue measure in $\Bbb{R}^n$.
\medskip

 Let $G$ and $G'$ be domains in $\Bbb{R}^n$, $n\geq 2$, and let
$Q:G\to[0,\infty]$ be a measurable function. A homeomorphism
$f:G\to{G'}$ is called a {\it $Q-$homeomorphism with respect to the
$p$-modulus} if \eqb \label{eq13.1}
 M_{p}(f\Gamma)\leq\int\limits_{G}
Q(x)\cdot\varrho^{p}(x)\ dm(x)\eqe for every family $\Gamma$ of
paths in $G$ and every admissible function $\varrho$ for $\Gamma$.

This conception is a natural generalization of the geometric
definition of a quasiconformal mapping: if $Q(x)\leq K<\infty$ a.e.,
then $f$ is quasiconformal under $p=n$, see 13.1 and 34.6 in
\cite{Va}, and quasiisometric under $1<p\neq n$, see \cite{Ge$_1$}.

 This class of $Q$-homeomorphisms with respect to the $n$-modulus was first considered
 in the papers \cite{MRSY$_1$}-\cite{MRSY$_3$}, see also the monograph \cite{MRSY}. The main goal of the theory of $
Q$-homeomorphisms is to clear up various interconnections between
properties of the majorant $Q(x)$ and the corresponding properties
of the mappings themselves. In particular, the problem of the local
and boundary behavior of $Q$-homeomorphisms has been studied in
$\Bbb{R}^n$ first in the case $Q\,\in\,BMO$ (bounded mean
oscillation) in the papers \cite{MRSY$_1$}-\cite{MRSY$_3$}  and then
in the case of $Q\,\in\,FMO$ (finite mean oscillation) and other
cases in the papers \cite{IR$_1$}, \cite{IR$_2$}, \cite{RS$_1$}.

\medskip

Note that the estimate of the type (\ref{eq13.1}) was first
established in the classical quasiconformal theory. Namely, it was
obtained in \cite{LV}, p. 221, for quasiconformal mappings in the
complex plane that \eqb\label{eq0.9} M(f\Gamma)\ \leq\
\int\limits_{\C}K(z)\cdot\rho^2(z)\ dxdy \eqe where
\eqb\label{eq0.10} K(z)\ =\
\frac{|f_z|+|f_{\overline{z}}|}{|f_z|-|f_{\overline{z}}|}\eqe is a
(local) maximal dilatation of the mapping $f$ at a point $z$.  Next,
it was obtain in \cite{BGMV}, Lemma 2.1, for quasiconformal mappings
in space, $n\geq2$, that \eqb\label{eq0.91}
 M(f\Gamma)\ \leq\ \int\limits_{D}K_I(x,f)\,\rho^n(x)\
dm(x) \eqe where $K_I(x,f)$ stands for the inner dilatation of $f$
at $x$, see (\ref{eq4.1.5}) below.

\medskip

Given a mapping $f:G \to \Rn$ with partial derivatives a.e.,
$f^\prime(x)$ denotes the Jacobian matrix of $f$ at $x \in D$ if it
exists, $J(x)=J(x,f)=\det f^\prime(x)$ the Jacobian of $f$ at $x$,
and $|f^\prime(x)|$ the operator norm of $f^\prime(x)$, i.e.,
$|f^\prime(x)|=\max \{ |f^\prime(x)h|: h \in \Rn, |h|=1\}.$ We also
let $l(f^\prime(x))= \min \{ |f^\prime(x)h|: h \in \Rn, |h|=1\}.$
The {\it outer dilatation}\index{outer dilatation $K_O(x,f)$} of $f$
at $x$ is defined by \eqb \label{eq4.1.4} K_O(x,f)= \left
\{\begin{array}{rr}
\frac{|f^\prime(x)|^n}{|J(x,f)|}, & {\rm if } \ J(x,f) \neq 0 \\
1, & {\rm if} \ f^\prime(x)=0 \\ \infty, & {\rm } \text{otherwise},
\end{array} \right. \eqe the {\it inner dilatation}\index{inner dilatation $K_I(x,f)$} of $f$ at $x$
by \eqb \label{eq4.1.5} K_I(x,f)= \left \{\begin{array}{rr}
\frac{|J(x,f)|}{l(f^\prime(x))^n}, &{\rm if } \ J(x,f) \neq 0 \\
1, & {\rm if } \ f^\prime(x)=0 \\
\infty, &{\rm} \text{otherwise},
\end{array} \right. \eqe

\medskip

The following notion generalizes and localizes the above notion of
$Q$--ho\-me\-o\-mor\-phism. It is motivated by the ring definition
of Gehring for qua\-si\-con\-for\-mal mappings, see, e.g.,
\cite{Ge$_2$}, introduced first by V. Ryazanov, U. Srebro, and E.
Yakubov in the plane and later on extended by V. Ryazanov and E.
Sevostyanov to the space case, see, e.g., \cite{RS$_2$}, \cite{RSY}
and Chapters 7 and 11 in \cite{MRSY}.

Let $E,$ $F\subset\Rn$ be arbitrary sets. Denote by $\Delta(E,F,G)$
a family of all curves $\gamma:[a,b]\rightarrow \overline{{\Bbb
R}^n}$ joining $E$ and $F$ in $G,$ i.e., $\gamma(a)\in E,\gamma(b)
\in F$ and $\gamma(t)\in G$ for $t \in (a, b).$ Given a domain $G$
in ${\Bbb R }^n,$ $n\ge 2,$ a (Lebesgue) measurable function
$Q:G\rightarrow\,[0,\infty]$, $x_0\in G,$
a homeomorphism $f:G\rightarrow \Rn$ is said to be a {\it ring
$Q$--homeomorphism at the point $x_0$} if
\eqb\label{eq1}
M_p\left(f\left(\Delta\left(S_1,\,S_2,\,A\right)\right)\right)\ \le
\int\limits_{A} Q(x)\cdot \eta^p(|x-x_0|)\ dm(x)
\eqe
for every ring $A=A(r_1,r_2,x_0)$ $=\{ x\,\in\,{\Bbb R}^n :
r_1<|x-x_0|<r_2\} $ and the spheres $S_i=S(x_0, r_i)=\{x\in {\Bbb
R}^n: |x-x_0|=r_i\}$, where $0<r_1<r_2< r_0\,\colon =\,{\rm dist}\,
(x_0,\partial D),$ and every measurable function $\eta :
(r_1,r_2)\rightarrow [0,\infty ]$ such that
\eqb\label{eq1n}\int\limits_{r_1}^{r_2}\eta(r)\ dr\ \ge\ 1\,.\eqe
$f$ is called a {\it ring $Q-$homeomorphisms with respect to the
$p$-modulus in the domain} $G$ if $f$ is a ring
$Q$--ho\-meo\-mor\-phism at every point $x_0\in G$.

Let $f:G\to \Rn, n\geq 2$ be  a quasiconformal mapping. The  {\it
angular dilatation} of the mapping $f:G\to \Rn$  at the point $x\in
G$ with  respect to $x_0\in G, x_0\neq x$ is defined by
\begin{equation}\label{eq1u*}
D_{f}(x,x_0)=\frac{J(x,f)}{l^n_f(x,x_0)}\,,
\end{equation}
where
$$
l_f(x,x_0)=\min\limits_{|h|=1} \frac{|\partial_h f(x)|}{|\langle
h,\frac{x-x_0}{|x-x_0|}\rangle|}.
$$
Here  $\partial_h f(x)$ denotes the derivative of $f$ at $x$ in the
direction  $h$ and the minimum is taken over all unit vectors $h\in
\Rn$, see \cite{GG}.

\medskip

We recall  that the estimate of the type (\ref{eq1}) was first
established in the classical quasiconformal theory in complex plane,
see \cite{GS}. Next, it was obtained in \cite{GG}, for
quasiconformal mappings in space, $n\geq2$, that \eqb\label{eq1d}
M\left(f\left(\Delta\left(S_1,\,S_2,\,A\right)\right)\right)\ \le
\int\limits_{A} D_f(x,x_0)\cdot \eta^n(|x-x_0|)\ dm(x)
\eqe
for every ring $A=A(x_0,r_1,r_2)$ $=\{ x\,\in\,{\Bbb R}^n :
r_1<|x-x_0|<r_2\} $ and the spheres $S_i=S(x_0, r_i)=\{x\in {\Bbb
R}^n: |x-x_0|=r_i\}$, where $0<r_1<r_2< r_0\,\colon =\,{\rm dist}\,
(x_0,\partial G),$ and every measurable function $\eta :
(r_1,r_2)\rightarrow [0,\infty ]$ such that (\ref{eq1n}) holds.

Note that, in particular, homeomorphisms $f:G\rightarrow \Rn$ in the
class $W_{loc}^{1,n}$ with $K_I(x,f)\in L_{loc}^1$ are ring
$Q$--ho\-me\-o\-mor\-phisms as well as $Q$--ho\-me\-o\-mor\-phisms
with $Q(x)=K_I(x,f),$ see, e.g., Theorem 6.10 and Corollary 4.9 in
\cite{MRSY$_1$}, or Theorem 4.1 in \cite{MRSY}.

\medskip

\large \cc
\section{Preliminaries}

Here a {\it condenser} is a pair $E\,=\,(A, C)$ where $A\subset
{\Bbb R}^n$ is open and $C$ is a non--empty compact set contained in
A\,. $E$ is a {\it ringlike condenser} if $B=A\setminus C$ is a
ring, i.e., if $B$ is a domain whose complement $\lRn\setminus B$
has exactly two components where $\lRn = \Rn\cup\{\infty\}$ is the
one point compactification of  $\Rn .$ $E$ is a {\it bounded
condenser} if $A$ is bounded. A condenser $E\,=\,(A, C)$ is said to
be in a domain $G$ if $A\subset G\,.$

The following proposition is immediate.

\medskip

\bpr \label{pr1}  If $f:G\rightarrow {\Bbb R}^n$ is open and
$E\,=\,(A, C)$ is a condenser in $G,$ then $\left(fA,\,fC\right)$ is
a condenser in $fG\,.$ \epr

\medskip

In the above situation we denote $fE\,=\,\left(fA,\,fC\right)\,.$
\bigskip

Let $E\,=\,\left(A,\,C\right)$ be a condenser. Then
$W_0(E)\,=\,W_0(A,\,C)$ denotes  the family of non--negative
functions $u:A\rightarrow R^1$ such that
(1)\,\,$u\,\in\,C_0(A),$\,\,(2)\,\, $u(x)\ge 1$ for $x\in C,$ and
(3)\,\,$u\,$ is $\,ACL\,.$  We set

\eqb\label{cap}
cap_p\,E\,\,=\,\,cap\,\left(A,\,C\right)\,\,=\,\,\inf
\limits_{u\,\in\,W_0\left(E\right) }\,\,\int\limits_{A}\,\vert\nabla
u\vert^p\,\,dm \eqe where
$$\vert\nabla
u\vert\,=\,{\left(\sum\limits_{i=1}^n\,{\left(\partial_i u\right)}^2
\right)}^{1/2}$$ and call the quantity (\ref{cap})  the {\it
$p$-capacity} of the condenser $E\,.$

\medskip

For the next statement,  see, e.g., \cite{Ge}, \cite{He} and
\cite{Sh}.

\bpr \label{pr2} Suppose  $E=(A,C)$ is a condenser such that $C$ is
connected. Then $${\rm cap}_p\,E=M_p(\Delta(\partial A,\partial C;
A\setminus C))\,.$$\epr

\medskip
We give here also the following two useful statements, see
Proposition 5 and 6 in \cite{Kr}.

\bpr \label{pr3} Let  $E=(A,C)$ be  a condenser such that $C$ is
connected. Then $$ {\rm cap}_p\,E\geqslant\frac{\left(\inf
m_{n-1}\,\sigma\right)^{p}}{\left[m(A\setminus C)\right]^{p-1}}$$
where  $ m_{n-1}\,\sigma$  denotes the $(n-1)$-dimensional area of
the  $C^{\infty}$-manifold $\sigma$ that is the boundary
$\sigma=\partial U$ of open set  $U$ containing $C$ and contained
along with its closure $\overline{U}$ in  $A$ and the infimum is
taken over all such $\sigma$.\epr

\medskip

\bpr \label{pr4} Let  $E=(A,C)$ be  a condenser such that $C$ is
connected. Then for $n-1<p\leq n$
$$\left(cap_{p}\,\,
E\right)^{n-1}\,\ge\,\gamma\,\frac{d(C)^{p}}{m(A)^{1-n+p}} $$ where
$\gamma$  is a positive constant that depends only on $n$ and $p$,
$d(A)$ is a diameter and $m(A)$ is the Lebesgue measure of $A$ in
$\Bbb{R}^n$. \epr

\medskip

\large \cc
\section{Characterization of ring $Q-$homeomorphisms with respect to
the $p$-modulus}

The theorems of this section extend the corresponding results in
\cite{RS$_2$}, see also Section 7.3 in the monograph \cite{MRSY},
from the case of $p=n$ to the case of $p\in (1,n]$. Below we use the
standard conventions: $a/\infty = 0$ for $a\ne\infty$ and
$a/0=\infty $ if $a>0$ and  $0\cdot\infty =0,$ see e.g. \cite{Sa},
p. 6.

\medskip

\blem{}\label{lem1} Let  $G$  be a domain in ${\Bbb R}^n\,, n\geq
2,$ $1<p\leq n$, $Q:G\rightarrow [0\,,\infty]\, $  a measurable
function and $q_{x_0}(r)$
 the mean of  $Q(x)$ over the sphere $|x-x_{0}|\,=\,r\,.$ Set
\eqb\label{eq2.3.7}
 I\ =\
I(x_0,r_1,r_2)\ =\ \int\limits_{r_1}^{r_2}\
\frac{dr}{r^{\frac{n-1}{p-1}}q_{x_0}^{\frac{1}{p-1}}(r)} \eqe and
$S_j=\{ x\in{\Bbb R}^n : |x-x_0|=r_j\},$ $j=1,2\,,$ where $x_0\in G$
and $0<r_1<r_2< r_0= dist \, (x_0,\partial G).$ Then whenever
$f:G\to\,{\Bbb R}^n$ is a ring $Q-$homeomorphism with respect to the
$p$-modulus at a point $x_0$ \eqb\label{eq2.3.8}
M_{p}\,\left(\Delta\left(fS_1,fS_2, fG\right)\right)\,\ \le\
\frac{\omega_{n-1}}{I^{p-1}} \eqe where $\omega_{n-1}$ is the area
of the unit sphere in ${\Bbb R}^n$. \elem

\medskip

{\it Proof.} With no loss of generality, we may assume that
$I\,\neq\,0$ because otherwise  $(\ref{eq2.3.8})$ is trivial, and
that $I\,\neq\,\infty$  because otherwise we can replace $Q(x)$ by
$Q(x)\,+\,\delta$  with arbitrarily small $\delta>0$ and then take
the limit as $\delta\to 0$ in $(\ref{eq2.3.8})$. The condition
$I\ne\infty$ implies, in particular, that  $q_{x_0}(r)\,\neq\,0$
a.e. in $(r_1\,,r_2)\,.$ Set
\eqb\label{eq2.3.9} \psi(t)= \left \{\begin{array}{rr}
1/[t^{\frac{n-1}{p-1}}q_{x_0}^{\frac{1}{p-1}}(t)]\ , & \ t\in
(r_1,r_2)\ ,
\\ 0\ ,  &  \ t\notin (r_1,r_2)\ .
\end{array} \right. \eqe
Then \eqb\label{eq2.3.10} \int\limits_{A}
Q(x)\cdot\psi^{p}(|x-x_0|)\ dm(x)\ =\ \omega_{n-1} I \eqe where $ A\
=\ A(x_0, r_1 ,r_2 )$.

Let  $\Gamma$  be a family of all paths joining the spheres  $S_1$
and $S_2$ in  $A\,.$ Let also  $\psi^{*}$ be a Borel function such
that $\psi^{*}(t)\,=\,\psi(t)$ for a.e. $t\,\in [0\,,\infty]\,.$
Such a function $\psi^{*}$ exists by the Lusin theorem, see, e.g.,
2.3.5 in \cite{Fu} and \cite{Sa}, p.
69. Then the function 
$$\rho(x)\ = \psi^*(|x-x_0|)/I$$
is admissible for $\Gamma_{}$ and since $f$ is a ring
$Q-$homeomorphisms with respect to the $p$-modulus we get by
$(\ref{eq2.3.10})$ that

\eqb\label{qe2.3.12} M_{p}(f\Gamma_{})\ \leq\ \int\limits_{A}
Q(x)\cdot {\rho_{}}^p (x) \ dm(x)\ =
\,\,\,\frac{\omega_{n-1}}{I^{p-1}}\, . \eqe
and the proof is complete.

\medskip

\medskip

The following lemma shows that the estimate $(\ref{eq2.3.8})$ cannot
be improved for ring $Q-$homeomorphisms with respect to the
$p$-modulus.

\medskip

\blem{}\label{lem2} Let $G$  be a domain in ${\Bbb R}^n$, $n\geq 2$,
$1<p\leq n$, $x_0\in G$, $0<r_1<r_2<r_0=dist(x_0,\partial G),$
$A=A(x_0,r_1,r_2)=\{ x\,\in\,{\Bbb R}^n : r_1<|x-x_0|<r_2\}$,
$Q:\,\to [0,\infty]$ be  a measurable function. Set
\eqb\label{eq2.3.13} \eta_0(r)=\frac{1}{Ir^{\frac{n-1}{p-1}}
q_{x_0}^{\frac{1}{p-1}}(r)} \eqe where $q_{x_0}(r)$ is the mean of
$Q(x)$ over the sphere $|x-x_0|=r$ and $I$ is given by
(\ref{eq2.3.7}). Then
\eqb\label{eq2.3.14} \frac{\omega_{n-1}}{I^{p-1}}=\int\limits_{A}
Q(x)\cdot \eta_0^{p}(|x-x_0|)\ dm(x)\le\int\limits_{A} Q(x)\cdot
\eta^{p}(|x-x_0|)\ dm(x) \eqe
whenever  $\eta :(r_1,r_2)\to [0,\infty]$ is  measurable and
\eqb\label{eq2.3.15} \int\limits_{r_1}^{r_2}\eta(r)\ dr\ = 1\, .
\eqe
\elem
\medskip

{\bf Proof.} If $I\,=\,\infty\,,$  then the left hand side in
$(\ref{eq2.3.14})$ is equal to zero and the inequality is obvious.
If $I\,=\,0\,,$ then
 $q_{x_0}(r)\,=\,\infty\,$ for a.e.
$r\,\in\,\left(r_1\,,r_2\right)$ and the both sides in
$(\ref{eq2.3.14})$ are equal to $\infty$. Hence we may assume below
that $0\,<\,I\,<\,\infty\,.$ Now, by $(\ref{eq2.3.13})$ and
$(\ref{eq2.3.15})$  $q_{x_0}(r)\neq 0$ and $\eta(r)\ne\infty$ a.e.
in $(r_1,r_2).$ Set
$\lambda(r)=r^{\frac{n-1}{p-1}}q_{x_0}^{\frac{1}{p-1}}(r)\eta(r)$
and $w(r)=1/r^{\frac{n-1}{p-1}}q_{x_0}^{\frac{1}{p-1}}(r)$. Then by
the standard conventions $\eta(r)=\lambda(r)w(r)$ a.e. in
$(r_1,r_2)$ and
\eqb\label{2.3.16}
 C\ \colon =\ \int\limits_{A} Q(x)\cdot \eta^p(|x-x_0|)\, dm(x)\,
=\, \omega_{n-1} \int\limits_{r_1}^{r_2}\lambda^{p}(r)\cdot w(r)\
dr\ . \eqe
By Jensen's inequality with weights, see, e.g., Theorem 2.6.2 in
\cite{Ra}, applied to the convex function $\varphi(t)\,=\,t^{p}$ in
the interval $\Omega\,=\,\left(r_1,r_2\right)$ with the probability
measure
\eqb\label{eq2.3.17}
 \nu \left(E\right)\ =\ \frac{1}{I}\ \int\limits_E w(r)\ dr
\eqe
we obtain that \eqb\label{eq2.3.18} \left(\dashint\lambda^p(r)w(r)\
dr\right)^{1/p}\ \ge\ \dashint\lambda(r)w(r)\ dr\ =\ \frac{1}{I}
\eqe
 where we also applied
that  $\eta(r)\,=\,\lambda(r)w(r)$ satisfies $(\ref{eq2.3.15})\,.$

Thus \eqb\label{eq2.3.19}
C\,\geq\,\frac{\omega_{n-1}}{I^{p-1}}
\eqe
and the proof of $(\ref{eq2.3.14})$ is complete.

\medskip

\newpage

Finally, combining Lemmas \ref{lem1} and \ref{lem2}, we obtain the
following statement.
\medskip

\bth{}\label{th1} Let $G$ be a domain in  ${\Bbb
R}^n\,,\,\,n\,\geq\,2\,,$ and $Q:G\,\rightarrow\,[0,\,\infty]$ a
measurable function. A homeomorphism $f:G\to {\Bbb R}^n$ is ring
$Q-$homeomorphism with respect to the $p$-modulus at a point
$x_0\,\in\,G$ if and only if for every $0<r_1<r_2< d_0=
dist\,(x_0,\partial G)$ \eqb\label{eq1.1.2aa}
M_{p}\left(\Delta\left(fS_1,\,fS_2\,
fG\right)\right)\,\leq\,\frac{\omega_{n-1}}{I^{p-1}}
\eqe
where  $S_1$ and $S_2$, $S_1=\{x\in \Bbb{R}^n: \ |x-x_0|=r_1\}$ and
$S_2=\{x\in \Bbb{R}^n: \ |x-x_0|=r_2\}$. $\omega_{n-1}$ is the area
of the unit sphere in ${\Bbb R}^n,\,\,$\,\,$I\ =\ I(x_0,r_1,r_2)\ =\
\int\limits_{r_1}^{r_2}\ \frac{dr}{r^{\frac{n-1}{p-1}}
q_{x_0}^{\frac{1}{p-1}}(r)}\,,$ $q_{x_0}(r)\,$ is the mean value of
$Q$ over the sphere $|x-x_0|\,=\,r\,.$ \eth

Note that  the infimum from the right hand side in $(\ref{eq1})$
holds for the function
\eqb\label{eq2.3.13a} \eta_0(r)=\frac{1}{Ir^{\frac{n-1}{p-1}}
q_{x_0}^{\frac{1}{p-1}}(r)}\,. \eqe

Theorem \ref{th1}  will have many applications in the theory of ring
$Q-$homeomorphisms with respect to the $p$-modulus, see, e.g., the
next section.


\large \cc
\section{ On finite  by Lipschitz mappimgs}
\medskip

Given a mapping $\varphi:E\to\Rn$ and a point $x\in E\subseteq\Rn$,
set  \eqb \label{eq8.1.6} L(x,\varphi)\ =\ \limsup_{y\to x\:y\in E}\
\frac{|\varphi(y)-\varphi(x)|}{|y-x|} \eqe and \eqb\label{eq8.1.7}
l(x,\varphi)\ =\ \liminf_{y\to x\:y\in E}\
\frac{|\varphi(y)-\varphi(x)|}{|y-x|}\eqe

Given a set $A\subseteq\Rn$, $n\geq1$, we say that a mapping
$f:A\to\Rn$ is called {\it Lipschitz}\index{Lipschitz mapping} if
there is number $L>0$ such that the inequality
\eqb\label{eq8.12.1}|\,f(x)-f(y)|\,\leq L\,|\,x-y|\eqe holds for all
$x$ and $y$ in $A$. Given an open set $\Omega\subseteq\Rn$, we say
that a mapping $f:\Omega\to\Rn$ is {\it finitely
Lipschitz}\index{finitely Lipschitz mapping} if $L(x,f)<\infty$  for
all $x\in\Omega$.

\medskip

\blem{} \label{lem3} Let $G$ and $G'$ be bounded domains in $\Rn$,
$n\geq 2$,  $Q:G\to [0,\infty]$ be a measurable function and  let
$f:G\to G'$ be a  ring $Q-$homeomorphism with respect to $p$-modulus
at a point $x_0\in G$, $1<p<n$. Then
\eqb\label{eqks*}m(fB(x_0,r_1))\leqslant \frac{c_{n,p}}{
I^{\frac{n(p-1)}{n-p}}(x_0,r_1,r_2)} \eqe for every $0<r_1<r_2< d_0=
dist\,(x_0,\partial G)$ where $I(x_0,r_1,r_2)$ is defined by
(\ref{eq2.3.7}) and $c_{n,p}$ is a positive constant that depends
only on $n$ and $p$. \elem

\medskip

{\bf Proof.} Let us consider  the condenser $(A_{t+\triangle
t},C_t)$, where $C_t=\overline{B(x_0,t)}, A_{t+\triangle
t}=B(x_0,t+\triangle t)$. Note that   $(fA_{t+\triangle t},fC_t)$ is
a  ringlike condenser in $\Bbb{R}^{n}$ and according to Proposition
\ref{pr2}, we have \eqb\label{eqks1.8}{\rm cap}_{p}(fA_{t+\triangle
t},fC_t)=M_{p}(\triangle(\partial fA_{t+\triangle t},\partial
fC_t;fR_{t})).\eqe In view of Proposition \ref{pr3}, we obtain
\eqb\label{eqks1.9}{\rm cap}_{p}\left(fA_{t+\triangle t},
fC_t\right)\geqslant\frac{\left(\inf
m_{n-1}\,\sigma\right)^{p}}{m\left(fA_{t+\Delta t}\setminus
fC_t\right)^{p-1}},\eqe where  $ m_{n-1}\,\sigma$  denotes the
$(n-1)$-dimensional area  of a   $C^{\infty}$-manifold $\sigma$ that
is the boundary  of an open set  $U$ containing $fC_t$ with its
closure $\overline{U}$ in $fA_{t+\triangle t}$  and the infimum is
taken  over all such $\sigma$.

On the other hand, by Lemma \ref{lem1}, we have  \eqb
\label{eqks1.10}M_{p}(\triangle(\partial fA_{t+\triangle t},\partial
fC_t;fR_{t}))\leqslant \frac{\omega_{n-1}}{\left(
\int\limits_{t}^{t+\Delta t}\frac{ds}{s^{\frac{n-1}{p-1}}\,
q_{x_0}^{\frac{1}{p-1}}(s)}\right)^{p-1}}\,. \eqe Combining
(\ref{eqks1.8})-(\ref{eqks1.10}), we obtain
$$\frac{\left(\inf m_{n-1}\, \sigma\right)^{p}}{m\left(fA_{t+\Delta
t}\setminus fC_t\right)^{p-1}}\leqslant \frac{\omega_{n-1}}{\left(
\int\limits_{t}^{t+\Delta t}\frac{ds}{s^{\frac{n-1}{p-1}}\,
q_{x_0}^{\frac{1}{p-1}}(s)}\right)^{p-1}}\,.$$

Applying the isoperimetric inequality to the numerator of the
fraction on the left-hand side  we came to the inequality

\eqb\label{eqks1.11}n\cdot\Omega_n^{\frac{1}{n}}
\left(m(fC_t)\right)^{\frac{n-1}{n}}\leq \omega^{\frac{1}{p}}_{n-1}
\left(\frac{m\left(fA_{t+\Delta t}\setminus
fC_t\right)}{\int\limits_{t}^{t+\Delta
t}\frac{ds}{s^{\frac{n-1}{p-1}}\,
q_{x_0}^{\frac{1}{p-1}}(s)}}\right)^{\frac{p-1}{p}}\eqe where
$\Omega_n$ is the volume of the unit ball in $\Rn$.

Now, setting $\Phi(t):=m\left(fB_t\right)$, we see from
(\ref{eqks1.11})  that

\eqb\label{eqks1.12}n\cdot\Omega_n^{\frac{1}{n}}\Phi^{\frac{n-1}{n}}(t)\leq
\omega^{\frac{1}{p}}_{n-1} \left(\frac{\frac{\Phi(t+\Delta
t)-\Phi(t)}{\Delta t}}{\frac{1}{\Delta t}\int\limits_{t}^{t+\Delta
t}\frac{ds}{s^{\frac{n-1}{p-1}}\,
q_{x_0}^{\frac{1}{p-1}}(s)}}\right)^{\frac{p-1}{p}} \,.\eqe  Since
the function $\Phi(t)$ is nondecreasing,  has finite derivative
$\Phi'(t)$ for a.e. $t$. Letting $\Delta t\to 0$ in (\ref{eqks1.12})
and taking into account that $\omega_{n-1}=n\Omega_n$, we obtain
\eqb\label{eqks1.14}\frac{n\Omega^{\frac{p-n}{n(p-1)}}_n}
{t^{\frac{n-1}{p-1}}\,q_{x_0}^{\frac{1}{p-1}}(t)}\leqslant
\frac{\Phi'(t)}{\Phi^{\frac{p(n-1)}{n(p-1)}}(t)}.\eqe

Integrating (\ref{eqks1.14}) under $1<p<n$ with respect to $t\in
[r_1,r_2]$, since
$$\int\limits_{r_1}^{r_2}\frac{\Phi'(t)}{\Phi^{\frac{p(n-1)}{n(p-1)}}(t)}\,dt\leqslant
\frac{n(p-1)}{p-n}\left(\Phi^{\frac{p-n}{n(p-1)}}(r_2)-\Phi^{\frac{p-n}{n(p-1)}}(r_1)\right),$$
see, e.g., Theorem  IV. 7.4 in \cite{Sa}, we observe that
\eqb\label{eqks1.15}\Omega^{\frac{p-n}{n(p-1)}}_n\int\limits_{r_1}^{r_2}\frac{dt}{t^{\frac{n-1}{p-1}}\,q_{x_0}^{\frac{1}{p-1}}(t)}\leqslant
\frac{p-1}{p-n}
\left(\Phi^{\frac{p-n}{n(p-1)}}(r_2)-\Phi^{\frac{p-n}{n(p-1)}}(r_1)\right).\eqe
From  (\ref{eqks1.15}) we conclude that
$$\Phi(r_1)\leqslant\left(\Phi^{\frac{p-n}{n(p-1)}}(r_2)+\Omega^{\frac{p-n}{n(p-1)}}_n\,
\frac{n-p}{p-1}\int\limits_{r_1}^{r_2}\frac{dt}{t^{\frac{n-1}{p-1}}\,q_{x_0}^{\frac{1}{p-1}}(t)}\right)^
{\frac{n(p-1)}{p-n}}$$ and hence

$$\Phi(r_1)\leqslant\Omega_n
\left(\frac{p-1}{n-p}\right)^{\frac{n(p-1)}{n-p}}\left(\int\limits_{r_1}^{r_2}\frac{dt}{t^{\frac{n-1}{p-1}}\,q_{x_0}^{\frac{1}{p-1}}(t)}\right)^
{-\frac{n(p-1)}{n-p}}\,.$$

\medskip

Combining  Lemmas \ref{lem2} and  \ref{lem3}, we have the following
statement.

\medskip

\blem{} \label{lem4}  Let $G$ and $G'$ be bounded domains in $\Rn$,
$n\geq 2$,  $Q:G\to [0,\infty]$ be a measurable function and let
$f:G\to G'$ be a ring $Q-$homeomorphism with respect to the
$p$-modulus. Then for $1<p<n$
\eqb\label{eqks*}m(fB(x_0,r_1))\leqslant c'_{n,p}\left[
\int\limits_{A(x_0,r_1,r_2)}\, Q(x)\eta^p(|x-x_0|)\, dm(x)
\right]^{\frac{n}{n-p}}\,.\eqe for every ring  $A=A(
x_0,r_1,r_2),$\,\, $0<r_1<r_2< d_0= dist\,(x_0,\partial G)$ and for
every measurable function  $\eta : (r_1,r_2)\to [0,\infty ]\,,$ such
that \eqb\label{2.13} \int\limits_{r_1}^{r_2}\eta(r)\ dr\ \geq\ 1\,.
\eqe where $c'_{n,p}$ is a positive constant that depends only on
$n$ and $p$. \elem

\medskip

\bth{} \label{th2}   Let $G$ and $G'$ be   domains in $\Rn$, $n\geq
2$, and  $Q:G\to [0,\infty]$ be  a measurable function such that

\eqb\label{3.1} Q_0=\overline{\lim\limits_{r \to 0}}\
\frac{1}{\Omega_{n}\varepsilon^{n}}\,\int\limits_{B(x_0,\varepsilon)
} Q(x)\, dm(x)<\infty. \eqe

Then for every   ring $Q-$homeomorphism  $f:G\to G'$  with respect
to the $p$-modulus,   $n-1<p<n$,  \eqb\label{eqks*}
L(x_0,f)=\limsup\limits_{x\to x_0}\frac{|f(x)-f(x_0)|}{|x-x_0|}\leq
\lambda_{n,p}\, Q^{\frac{1}{n-p}}_{0}\eqe where $\lambda_{n,p}$ is a
positive constant that depends only on  $n$ and $p$. \eth

\medskip

{\it Proof.} Let us consider  the spherical ring $A(x_0,\varepsilon,
2\varepsilon)=\{x:\ \varepsilon<|x-x_0|< 2\varepsilon\}$, $x\in G,$
ñ $\varepsilon>0$ such that $A(x_0,\varepsilon, 2\varepsilon)\subset
G$. Since
$\left(fB\left(x_0,2\varepsilon\right),\overline{fB\left(x_0,\varepsilon\right)}\right)=\left(fB\left(x_0,2\varepsilon\right),f\overline{B\left(x_0,\varepsilon\right)}\right)
$ are ringlike condensers in  $G'$ and, according to Proposition
\ref{pr2}, we obtain
$$cap_{p}\ (fB(x_0,2\varepsilon),\overline{fB(x_0
,\varepsilon)})=M_{p}(\triangle(\partial
fB(x_0,2\varepsilon),\partial f
B(x_0,\varepsilon);fA(x_0,\varepsilon, 2\varepsilon)))\,.$$ Note
that, in view of the homeomorphism of $f,$
$$\triangle\left(\partial
fB\left(x_0,2\varepsilon\right),\partial
fB\left(x_0,\varepsilon\right);fA(x_0,\varepsilon,
2\varepsilon)\right)=f\left(\triangle\left(\partial
B(x_0,2\varepsilon) ,\partial B(x_0,\varepsilon);A(x_0,\varepsilon,
2\varepsilon)\right)\right).$$

By Proposition \ref{pr4}  \eqb\label{2.5} cap_{p}\
(fB(x_0,2\varepsilon),\overline{fB(x_0,\varepsilon)}) \geq
\left(\gamma\frac{d^{p}(fB(x_0,\varepsilon))}{m^{1-n+p}(fB(x_0,2\varepsilon))}\right)^{\frac{1}{n-1}}
\eqe where    $\gamma$  is a positive constant that depends only on
$n$ and $p$,  $d(A)$ is the  diameter and $m(A)$ is the Lebesgue
measure of $A$ in $\Bbb{R}^n$.

By  the definition of ring $Q-$homeomorphisms with respect to the
$p$-modulus \eqb\label{2.5eq} cap_{p}\
(fB(x_0,2\varepsilon),\overline{fB(x_0,\varepsilon)}) \leq
\frac{1}{\varepsilon^p}\int\limits_{A(x_0,\varepsilon,
2\varepsilon)} Q(x)\ d\,m(x) \eqe because the function
$$ \eta(t)\,=\,\left
\{\begin{array}{rr} \frac{1}{\varepsilon}, & {\rm if } \ t\in (\varepsilon,2\varepsilon), \\
0, & {\rm if} \ t\in \Bbb{R}\setminus (\varepsilon,2\varepsilon)
\end{array}\right.
$$
satisfies (\ref{eq1n}) for $r_1=\varepsilon$ and $r_2=2\varepsilon$.

Next, the  function

$$\widetilde{\eta}(t)\,=\,\left
\{\begin{array}{rr} \frac{1}{2\varepsilon}, & {\rm if } \ t\in (2\varepsilon,4\varepsilon) \\
0, & {\rm if} \ t\in \Bbb{R}\setminus (2\varepsilon,4\varepsilon),
\end{array}\right.
$$
satisfies (\ref{eq1n}) for $r_1=2\varepsilon$ and $r_2=4\varepsilon$
and hence by Lemma \ref{lem4} we have the following estimates:

\eqb\label{eqks*}m(fB(x_0,2\varepsilon))\leqslant
c''_{n,p}\,\varepsilon^n\left[
\frac{1}{m(B(x,4\varepsilon))}\int\limits_{B(x_0,4\varepsilon)}\,
Q(x)\, dm(x) \right]^{\frac{n}{n-p}}\,,\eqe where $c''_{n,p}$ is a
positive constant that depends only on $n$ and $p$.

Combining  (\ref{eqks*}), (\ref{2.5eq}) and  (\ref{2.5}), we obtain
$$\frac{d(fB(x_0,\varepsilon))}{\varepsilon}\leq\lambda_{n,p}\left(\frac{1}{m(B(x_0,4\varepsilon))}
\int\limits_{B(x_0,4\varepsilon)} Q(y)\,
dy\right)^{\frac{n(1-n+p)}{p(n-p)}}\left[
\frac{1}{m(B(x,2\varepsilon))}\int\limits_{B(x_0,2\varepsilon)}\,
Q(x)\, dm(x) \right]^{\frac{n-1}{p}}$$ and hence
$$
L(x_0,f)=\limsup\limits_{x\to
x_0}\frac{|f(x)-f(x_0)|}{|x-x_0|}\leq\limsup\limits_{\varepsilon\to
0}\frac{d(fB(x_0,\varepsilon))}{\varepsilon}\leq \lambda_{n,p}\,
Q^{\frac{1}{n-p}}_{0} ,
$$
where $\lambda_{n,p}$ is a positive constant that depends only on
$n$ and $p$.

\bigskip

\bcor \label{cor1}  Let $G$ and $G'$ be domains in $\Rn$, $n\geq 2$,
$f:G\to G'$ be a ring $Q-$homeomorphism with respect to the
$p$-modulus, $n-1<p<n$, such that

\eqb\label{3.1} \overline{\lim\limits_{r \to 0}}\
\frac{1}{\Omega_{n}\varepsilon^{n}}\,\int\limits_{B(x_0,\varepsilon)
} Q(x)\, dm(x)<\infty \ \ \forall \ x_0\in G. \eqe Then  $f$ is
finitely Lipschitz. \ecor

\medskip

Note that the theory of ring $Q-$homeomorphisms with respect to
$p$-modulus can be applied to mappings in the Orlich-Sobolev classes
$W_{loc}^{1,\varphi}$ with a Calderon type condition on $\varphi$
and, in particular, to the Sobolev classes $W_{loc}^{1,p}$ for
$p>n-1$,  cf. \cite {KRSS}.

\bigskip
\bf{ ACKNOWLEDGMENT:} I thank Professor Vladimir Ryazanov for
interesting discussions and valuable comments.


\medskip
\noindent Ruslan Salimov, \\
Institute of Applied Mathematics and Mechanics,\\
National Academy of Sciences of Ukraine, \\
74 Roze Luxemburg str., 83114 Donetsk, UKRAINE \\
Phone: +38 -- (062) -- 3110145 Fax: \ \ \ \,+38 -- (062) -- 3110285 \\
salimov07@rambler.ru, ruslan623@yandex.ru

\end{document}